\documentclass[11pt,a4paper,fleqn]{article}
\usepackage{amsmath}
\usepackage{amssymb}
\usepackage{amscd}
\usepackage[dvips]{graphicx}

\newtheorem{theorem}{Theorem}

\newtheorem{corollary}[theorem]{Corollary}
\newtheorem{definition}[theorem]{Definition}
\newtheorem{example}[theorem]{Example}
\newtheorem{lemma}[theorem]{Lemma}
\newtheorem{proposition}[theorem]{Proposition}
\newtheorem{remark}[theorem]{Remark}
\newenvironment{proof}{{\bf Proof. }}{\hfill$\rule{1ex}{1ex}$\par\medskip}
\begin{document}
\newcommand{\bt}{\begin{theorem}}
\newcommand{\et}{\end{theorem}}
\newcommand{\bd}{\begin{definition}}
\newcommand{\ed}{\end{definition}}
\newcommand{\bs}{\begin{proposition}}
\newcommand{\es}{\end{proposition}}
\newcommand{\bp}{\begin{proof}}
\newcommand{\ep}{\end{proof}}
\newcommand{\be}{\begin{equation}}
\newcommand{\ee}{\end{equation}}
\newcommand{\ul}{\underline}
\newcommand{\br}{\begin{remark}}
\newcommand{\er}{\end{remark}}
\newcommand{\bex}{\begin{example}}
\newcommand{\eex}{\end{example}}
\newcommand{\bc}{\begin{corollary}}
\newcommand{\ec}{\end{corollary}}
\newcommand{\bl}{\begin{lemma}}
\newcommand{\el}{\end{lemma}}
\title{Bruck decomposition for endomorphisms of quasigroups}
\author{P\'eter T. Nagy (Debrecen) and Peter Plaumann (Cuernavaca)}
\footnotetext{2000 {\em Mathematics Subject Classification:\/20N05} .}
\footnotetext{{\em Key words and phrases:} Bruck extension of quasigroups, non-associative extension of groups, special varieties of quasigroups, {\em LF}-quasigroups\par
This research was supported by the Hungarian Scientific Research Fund (OTKA), Grant K 67617  
and by the Program {\em C\'atedras Especiales. SRE} of the Mexican Government.}
\maketitle
\begin{abstract}
In the year 1944 R. H. Bruck has described a very general construction method which he called 
the extension of a set by a quasigroup. 
We use it to construct a class of examples for {\em LF}-quasigroups in which the image of the map 
$e(x)=x\backslash x$ is a group.  More generally, we consider the variety of quasigroups which 
is defined by the property that the map $e$ is an endomorphism and its subvariety 
where the image of the map $e$ is a group. We characterize quasigroups belonging to 
these varieties using their Bruck decomposition with respect to the map $e$. 
\end{abstract}

\section{Introduction}

A binary algebra $(Q,\cdot)$ with multiplication $(x,y)\mapsto x\cdot y$ 
is called a {\it quasigroup} if the equations $a\cdot y = b$ and
$x\cdot a = b$ have precisely one solution in $Q$ which we denote
by $y =a\backslash b$ and $x = b/a$. The elements  $1_{l}(a)=a/a$
respectively $1_{r}(a)=a\backslash a$ are the left respectively the right local
unit element of the element $a$. If the left (right) local
unit elements coincide for all elements of $(Q,\cdot)$,
then the element $1_{l}=1_{l}(a)$, resp. $1_{r}=1_{r}(a)$, is called the left, resp.  
right, unit element of $(Q,\cdot)$. If a quasigroup $(Q,\cdot)$ has both left 
and right unit elements then they coincide $1=1_{l}=1_{r}$, in this case 
$(Q,\cdot)$ is called  a {\it loop}.\\[1ex]\indent 
In the year 1944 R. H. Bruck has described a very general construction method which he called 
the extension of a set by a quasigroup (cf. \cite{Bru1}, \cite{Bru2}). 
Epimorphisms of quasigroups in general cannot be described by cosets of a 
normal subquasigroup, but only by congruence relations in the sense of 
universal algebra. Bruck's construction takes this into account giving a manageable 
description of quasigroup epimorphisms. In this note we discuss this method for 
endomorphisms of quasigroups.  \\[1ex]\indent
A quasigroup $(Q,\cdot)$ is called an {\em LF-quasigroup} if the identity 
\[x\cdot yz = xy\cdot (x\backslash x\cdot z)\]
holds in $Q$. In his book \cite {Bel1} V. D. Belousov initiated a systematic study 
of {\em LF}-quasigroups using isotopisms. Recently progress has been made in this topic 
(cf. \cite{Shc}, \cite{KKP}, \cite{PSS}). It is known that in an {\em LF}-quasigroup 
the map $e(x)=x\backslash x$ is an endomorphism, which we call the left deviation. 
In this situation Bruck's theory is available. 
We use it to construct a class of examples for {\em LF}-quasigroups $Q$ in which $e(Q)$ is 
a group.  \\[1ex]\indent
More generally, we consider the variety ${\mathfrak D}_l$ of quasigroups which 
is defined by the property that the left deviation is an endomorphism and its subvariety ${\mathfrak aD}_l$ 
where the image of the left deviation is a group. We characterize quasigroups belonging to 
these varieties using their Bruck decomposition with respect to their left deviation. 

\section{The Bruck decomposition of a quasigroup}

In Bruck's papers \cite{Bru1}, (Theorem 10\,A, pp. 166-168) and \cite{Bru2}, (pp. 778-779) 
a  principal construction for quasigroups is given. Let 
$(E,\cdot,\backslash, /)$ be a quasi\-group, let $T$ be a set and let 
$\{\triangledown_{a,b};\; a,b\in E\}$ be a family of multiplications on $T$.  
Define on the set $Q=T\times E$ a multiplication by 
\begin {equation}\label{mult}(\alpha,a)\circ(\beta,b)=
(\alpha\,\triangledown_{a,b}\,\beta,ab), \quad (\alpha,a),(\beta,b)\in T\times E.\end{equation}
Then $(Q,\circ)$ is a quasigroup if and only if for any $a,b\in E$ the multiplication 
$\triangledown_{a,b}$ on $T$ defines a quasigroup $\mathcal{T}_{a,b}=(T,\triangledown_{a,b})$. 
In this case we call ${\mathcal B}=\big(E,T,(\triangledown_{a,b})_{a,b\in E}\big)$ 
{\em Bruck system} and put $Q({\mathcal B})=Q$. Obviously, the projection 
$((\alpha,a)\mapsto a):Q({\mathcal B})\to E$ is an epimorphism of quasigroups. 
(Cf. \cite{Che}, pp. 35-36.) We call this epimorphism the {\em canonical epimorphism} for ${\mathcal B}$.\\[1ex]\indent
Conversely, let $(Q,\circ)$ and $(E,\cdot)$ be quasigroups and let $\pi:Q\to E$ be an epimorphism. 
For the inverse images $T_a=\pi^{\leftarrow}(a), \;a\in E,$ one has 
\[T_a\circ T_b = T_{ab},\]
hence the set $\{T_a, \;a\in E\}$ forms a quasigroup $\mathcal{E'}$ isomorphic to $(E,\cdot)$. 
Using a transversal for the partition $\{T_a, \;a\in E\}$ of the set $Q$ one can identify the 
inverse images $T_a, \;a\in E,$ with a subset $T\subset Q$ and the set $Q$ with the cartesian product 
$T\times E$.  The multiplication in $Q$ can be written in the form (\ref{mult}), where 
$\alpha, \beta\in T,\; a,b\in E$. This means that we have obtained a Bruck system 
${\mathcal B}_{\pi}=\big(E,T,(\triangledown_{a,b})_{a,b\in E}\big)$ for which 
$Q({\mathcal B}_{\pi})$ is isomorphic to $Q$ and which has $\pi$ as canonical epimorphism. 
We call this representation of $(Q,\circ)$ a {\em Bruck 
decomposition of $(Q,\circ)$ with respect to $\pi:Q\to E$}. \bigskip

Let $\eta :Q\to Q$ be an endomorphism. We consider the set $E$ of the congruence classes 
$T_x=\eta^{\leftarrow}(x), \;x\in \eta(Q)$. The multiplication $T_x\star T_y = T_{xy}$ 
defines a quasigroup $(E,\star)$ such that the mapping $\pi :(x\mapsto T_x):Q\to E$ is an 
epimorhism. We put $\iota = \eta\circ\pi^{\leftarrow}=(T_x\mapsto\eta(x)):E\to Q$. Then 
$\iota$ is an injective homomorphism (see Theorem 6.12, \cite{BS}, p. 50). \\
We consider the Bruck decomposition $Q=T\times E$ with respect to the epimorphism $\pi$.  
For $a\in E=\{T_x, \;x\in \eta(Q)\}$ one has $\iota(a)=(\gamma(a),g(a))$, where 
$\gamma(a)\in T,\; g(a)\in E$. The maps $\gamma: E\to T$ and $g: E\to E$ satisfy 
\[(\gamma(ab),g(ab))=\iota(ab)=\iota(a)\iota(b)=\]
\[=(\gamma(a),g(a))(\gamma(b),g(b))=(\gamma(a)\,\triangledown_{g(a),g(b)}\,\gamma(b),g(a)g(b)).\]
It follows that $g: E\to E$ is an endomorphism and 
\begin {equation}\label{hom}\gamma(ab)=\gamma(a)\,\triangledown_{g(a),g(b)}\,\gamma(b)\end{equation}
holds for all $a,b\in E$. We call the structure described here the {\em Bruck decomposition 
of the quasigroup $(Q,\circ)$ with respect to the endomorphism $\eta:Q\to Q$}. \\[1ex]\indent 
If the quasigroup $Q$ is a loop and $\eta :Q\to Q$ is an endomorphism, then $K=\eta^{\leftarrow}(1)$ 
is a normal subloop of $Q$. In this situation $Q$ is a semidirect product of $K$ and $\eta(Q)$ 
if and only if $\eta = \eta ^2$ holds. We describe the Bruck decomposition 
with respect to an idempotent endomorphism for arbitrary quasigroups. 
\bs An endomorphism $\eta$ of a quasigroup $Q$ is idempotent if and only if in the Bruck 
decomposition with respect to $\eta $ the maps $\gamma: E\to T$ and $g: E\to E$ 
satisfy \[g^2=g\quad \text{and}\quad \gamma\circ g=\gamma, \]
i. e. if and only if the endomorphism $g: E\to E$ is idempotent and the map $\gamma: E\to T$ 
factors over the congruence relation defined by $g$ on $E$. 
\es
\bp Since $\eta(\alpha,a)=\iota(\pi(\alpha,a))=\iota(a)=(\gamma(a),g(a))$ the assertion follows from 
\[\eta(\eta(\alpha,a))=\eta(\gamma(a),g(a))=(\gamma(g(a),g(g(a))).\] 
\ep

\section{The left deviation of a quasigroup} 
For a quasigroup $(Q,\cdot,\backslash , /)$ we call the map $e=(x\mapsto x\backslash x):Q\to Q$ the 
{\em left deviation}. As mentioned in the preliminaries the left deviation of $x\in Q$ is the local 
right unit element of $x$. In a Bruck decomposition $Q=T\times E$ with respect to an epimorphism 
$Q\to E$ the deviation is $e(\alpha,a)=(\alpha\backslash\alpha ,a\backslash a)$, where  
$a\backslash a$ is computed in the quasigroup $E$ and $\alpha\backslash\alpha$ is computed in the 
quasigroup $\mathcal{T}_{a,a\backslash a}=(T,\triangledown_{a,a\backslash a})$. \\
Obviously, the quasigroups in which {\em the left deviation is an endomorphism 
form a variety} ${\mathfrak D}_l$ of quasigroups. For $(Q,\cdot,\backslash , /)\in {\mathfrak D}_l$ 
we consider the Bruck decomposition $Q=T\times E$ with respect to the left deviation. 
In this case $e(\alpha,a)=(\gamma(a),g(a))$ and hence $g(a)=a\backslash a$ (computed in $E$) 
and \begin{equation}\label{idem}\alpha\,\triangledown_{a,g(a)}\,\gamma(a)=\alpha.\end{equation} 
\bt \label{Dl} A quasigroup $Q$ belongs to the variety ${\mathfrak D}_l$ if and only if 
there exists a Bruck system ${\mathcal B}=\big(E,T,(\triangledown_{a,b})_{a,b\in E}\big)$  
satisfying 
\begin{enumerate} 
\item[{\em (i)}] $Q\cong Q({\mathcal B})$, 
\item[{\em (ii)}]  the quasigroup $\mathcal{T}_{a,a\backslash a}=(T,\triangledown_{a,a\backslash a})$ 
has a right unit element, denoted by $\epsilon (a)$, for any $a\in E$,
\item[{\em (iii)}] the map $(a\mapsto(\epsilon (a),a\backslash a)):E\to Q({\mathcal B})$ is a homomorphism. 
\end{enumerate}
In this case the left deviation of $Q({\mathcal B})$ is the map \[e=\left((\alpha,a)\mapsto(\epsilon (a),a\backslash a)\right):Q({\mathcal B})\to Q({\mathcal B}).\]
\et 
\bp Assume first that $Q$ belongs to ${\mathfrak D}_l$ and consider the Bruck decomposition 
with respect to the left deviation $e(x)=x\backslash x$ of $Q$. Then $E$ is isomorphic to the 
subquasigroup $e(Q)$ of $Q$. Putting $\epsilon = \gamma$ the assertion (ii) 
follows from equation (\ref{idem}) and the assertion (iii) follows from equation (\ref{hom}).\\
Conversely, if ${\mathcal B}=\big(E,T,(\triangledown_{a,b})_{a,b\in E}\big)$ is a Bruck system 
satisfying (i) 
then $e(\alpha ,a)=(\alpha ,a)\backslash(\alpha ,a)=(\alpha ',a\backslash a)$, where 
$\alpha = \alpha \,\triangledown_{a,a\backslash a}\,\alpha '$. From (ii) it follows that 
$\alpha '= \epsilon(a)$ and the deviation satisfies $e(\alpha,a)=(\epsilon (a),a\backslash a)$. 
Hence we obtain from (iii) that the quasigroup $Q$ belongs to the variety ${\mathfrak D}_l$. \ep
\bc Let ${\mathcal B}=\big(E,T,(\triangledown_{a,b})_{a,b\in E}\big)$ be a Bruck system 
satisfying the conditions {\em (ii)} and {\em (iii)} of the previous theorem. Then ${\mathcal B}$ 
is a Bruck decomposition of the quasigroup $Q({\mathcal B})$ with respect to the left deviation  
of $Q({\mathcal B})$ if and only if the homomorphism 
\[\iota =(a\mapsto(\epsilon (a),a\backslash a)):E\to Q({\mathcal B})\] is injective.   
\ec
\bex \label{ex1} {\em Let $(E,\cdot)$, $\mathcal{T}^{(1)}=(T,\circ)$ and $\mathcal{T}^{(2)}=(T,\star)$ be 
quasigroups such that the following properties are satisfied:}
\begin{enumerate}    
\item[{\em (i)}] $E$ is a ${\mathfrak D}_l$-quasigroup,
\item[{\em (ii)}] $\mathcal{T}^{(1)}$ has a right unit element $\epsilon $,  
\item[{\em (iii)}] $\epsilon $ is idempotent in the quasigroup $\mathcal{T}^{(2)}$.
\end{enumerate}
{\em Put} 
\begin{enumerate}    
\item[{\em (i)}] $\mathcal{T}_{a,a\backslash a}=(T,\triangledown_{a,a\backslash a})$=$\mathcal{T}^{(1)}=(T,\circ)$,
\item[{\em (ii)}] $\mathcal{T}_{a\backslash a,b\backslash b}=(T,\triangledown_{a\backslash a,b\backslash b})$=
$\mathcal{T}^{(2)}=(T,\star)$ if $b\backslash b\neq (a\backslash a)\backslash(a\backslash a)$,  
\item[{\em (iii)}] $\mathcal{T}_{c,d}=(T,\triangledown_{c,b\backslash d})$ arbitrary in all other cases.
\end{enumerate}
{\em According to Theorem \ref{Dl} the multiplication $(\alpha,a)\circ(\beta,b)=
(\alpha\,\triangledown_{a,b}\,\beta,ab)$ on the set $T\times E$ is a ${\mathfrak D}_l$-quasigroup $Q$ with 
left deviation $e(\alpha,a)=(\epsilon ,a\backslash a)$. For $(E,\cdot)$
one can take idempotent quasigroups or groups. Clearly, the construction of the quasigroup $Q=T\times E$ gives a 
Bruck decomposition with respect to the left deviation of $Q$ if and only if the left deviation of $(E,\cdot)$ 
is an automorphism. If $(E,\cdot)$ is a group with unit element $1$ then the left deviation of $Q$ is the 
constant $(\epsilon ,1)$ which is the right unit element of $Q$.} 
\eex

\section{Associative image of the left deviation map} 

The class of ${\mathfrak D}_l$-quasigroups for which {\em the image of the left deviation map is a group} forms a variety, too, 
as can be seen from the identities:  
\begin{equation}
xy\backslash xy = x\backslash x \cdot y\backslash y, \quad (x\backslash x \cdot y\backslash y)\cdot z\backslash z = x\backslash x \cdot (y\backslash y\cdot z\backslash z).
\end{equation}
We denote this variety by ${\mathfrak {aD}}_l$ and investigate their Bruck decomposition. As an immediate consequence of 
Theorem \ref{Dl} one obtains 
\bt \label{adl}A quasigroup $Q$ belongs to the variety ${\mathfrak aD}_l$ if and only if 
there exists a Bruck system ${\mathcal B}=\big(E,T,(\triangledown_{a,b})_{a,b\in E}\big)$  
with $Q\cong Q({\mathcal B})$ satisfying 
\begin{enumerate}    
\item[{\em (i)}] $E$ is a group (with unit element $1$),
\item[{\em (ii)}] the quasigroup $\mathcal{T}_{a,1}=(T,\triangledown_{a,1})$ 
has right unit element, denoted by $\epsilon (a)$, for any $a\in E$,
\item[{\em (iii)}] the map $\epsilon :E\to \pi^{\leftarrow}(1)=\mathcal{T}_{1,1}=(T,\triangledown_{1,1})$ is a homomorphism of the group $E$ into the normal subquasigroup $\pi^{\leftarrow}(1)=\mathcal{T}_{1,1}=(T,\triangledown_{1,1})$ of $Q({\mathcal B})$.
\end{enumerate}
In this case one has $e(\alpha,a)=(\epsilon (a),1)$ for any $(\alpha ,a)\in T\times E$. \hfill$\rule{1ex}{1ex}$
\et  
\bp $E$ is a group since it is isomorphic to the image of the left deviation map. Hence  
$e(\alpha,a)=(\epsilon (a),a\backslash a)=(\epsilon (a),1)$ for any $(\alpha ,a)\in T\times E$ and  
$\epsilon :(a\mapsto(\epsilon (a),1)):E\to \pi^{\leftarrow}(1)=\mathcal{T}_{1,1}$. 
\ep
\bex \label{ex2} {\em Let $(E,\cdot)$, $(T,\circ)$ be quasigroups and let 
$\epsilon :E\to T$  be a homomorphism such that the following properties are satisfied:}
\begin{enumerate}    
\item[{\em (i)}] $(E,\cdot)$ is a group with unit element $1$,
\item[{\em (ii)}] $\epsilon(1)$ is a right unit element of $(T,\circ)$.
\end{enumerate}
{\em Put} 
\begin{enumerate}    
\item[{\em (a)}] $\alpha\,\triangledown_{a,1}\,\beta = \left(\alpha /\epsilon(a)\right)\circ\beta$ 
in the quasigroup $\mathcal{T}_{a,1}=(T,\triangledown_{a,1})$, where $/$ is the right division in $(T,\circ)$,
\item[{\em (b)}] $\mathcal{T}_{a,b}=(T,\triangledown_{a,b})$ arbitrary for $b\neq 1$.
\end{enumerate}
{\em According to Theorem \ref{adl} the multiplication $(\alpha,a)\circ(\beta,b)=
(\alpha\,\triangledown_{a,b}\,\beta,ab)$ on the set $T\times E$ is an ${\mathfrak {aD}}_l$-quasigroup. 
For $(E,\cdot)$ and $(T,\circ)$ one can take groups having $\epsilon :E\to T$ as a group homomorphism.
The decomposition $Q=T\times E$ is a Bruck decomposition with respect to the left deviation of $Q$ 
if and only if the homomorphism $\epsilon :E\to T$ is injective.} 
\eex

\section{{\em LF}-quasigroups}

It is known that the {\em LF}-quasigroups form a subvariety of ${\mathfrak D}_l$ (cf. \cite{Bel2} p. 108. and Lemma 2.1. in \cite{Shc}). We shall now give examples of {\em LF}-quasigroups even belonging to the variety ${\mathfrak aD}_l$. 
\bex \label{ex3} {\em Let $E$, $T$ be groups and let 
$\epsilon :E\to T$  be a homomorphism.} 
{\em Put \[\alpha\,\triangledown_{a,b}\,\beta = \alpha \cdot\epsilon(a)^{-1}\cdot\beta\] 
for all $a,b\in E$. Then every $\mathcal{T}_{a,b}=(T,\triangledown_{a,b})$ is a group 
(with unit element $\epsilon(a)$) which is isotopic and hence isomorphic to the group $T$.}\\
{\em As in the previous example the multiplication $(\alpha,a)\circ(\beta,b)=
(\alpha\,\triangledown_{a,b}\,\beta,ab)$ on the set $T\times E$ defines an ${\mathfrak {aD}}_l$-quasigroup 
(with left unit element $(1,1)$) in which the left deviation is given by $e(\alpha,a)=(\epsilon(a),1)$  (Theorem \ref{adl}).
\bt Let $E$, $T$ be groups and let $\epsilon :E\to T$  be a homomorphism. The set $Q=T\times E$ 
equipped with multiplication $(\alpha,a)\circ(\beta,b)=(\alpha \cdot\epsilon(a)^{-1}\cdot\beta,ab)$ 
is an LF-quasigroup with left unit element $(1,1)$ satisfying the left inverse property. 
\et
\bp An easy calculation shows 
\[(\alpha,a)\circ\big((\beta,b)\circ(\gamma,c)\big) = 
(\alpha\cdot\epsilon(a)^{-1}\cdot\beta\cdot\epsilon(b)^{-1}\cdot\gamma,abc).\]
On the other hand 
\[\big((\alpha,a)\circ(\beta,b)\big)\circ\big((\epsilon(a),1)\circ(\gamma,c)\big) = 
(\alpha\cdot\epsilon(a)^{-1}\cdot\beta ,ab)\circ(\epsilon(a)\cdot\gamma,c)=\]
\[=(\alpha\cdot\epsilon(a)^{-1}\cdot\beta\cdot\epsilon(ab)^{-1}\cdot\epsilon(a)\cdot\gamma,abc).\]
Hence $(Q,\circ)$ is an {\em LF}-quasigroup in which $(1,1)$ is the left unit element. The 
left inverse of an element $(\alpha ,a)$ is given by 
$(\epsilon(a)\cdot{\alpha}^{-1}\cdot\epsilon(a)^{-1},a^{-1})$. Indeed one has 
$(\epsilon(a)\cdot{\alpha}^{-1}\cdot\epsilon(a)^{-1},a^{-1})\circ\big((\alpha ,a)\circ(\beta ,b)\big) = 
(\beta ,b)$.
\ep
The decomposition $Q=T\times E$ is a Bruck decomposition with respect to the left deviation of $Q$ 
if and only if the homomorphism $\epsilon :E\to T$ is injective. }
\eex
We note that in accordance with Theorem 4.1 in \cite{Shc} the quasigroup $Q=T\times E$ is isotopic to the 
direct product of the groups $T$ and $E$, the isotopism is given by the triple $(\phi ,\text{id,id})$, 
where $\phi :(\alpha ,a)\mapsto (\alpha\cdot\epsilon(a)^{-1},a)$.

P\' eter T. Nagy\\
Institute of Mathematics\\
University of Debrecen\\
P.O.B. 12\\
4010 DEBRECEN, HUNGARY \vspace{1mm}\\
{\it E-mail}: {\tt {}nagypeti@math.unideb.hu}\vspace{4mm}\\
Peter Plaumann\\
Mathematisches Institut der\\
Universit\"at Erlangen-N\"urnberg\\
Bismarckstr.~$1 \frac 12$\\
91054 ERLANGEN, GERMANY \vspace{1mm}\\
{\it E-mail}: {\tt {}peter.plaumann@mi.uni-erlangen.de}
\end{document}